\pdfoutput=1
\documentclass[12pt,a4paper]{style}

\pagespan{1}{?}

\newtheorem{theorem}{Theorem}[section]

\newtheorem{corollary}[theorem]{Corollary}

\theoremstyle{remark}

\usepackage{times}
\usepackage{bbding}

\usepackage{hyperref}
\usepackage{float}
\usepackage{caption}
\usepackage{subfig}

\usepackage{tikz}
\usetikzlibrary{shapes,arrows}

\newcommand{\lisays}{\vspace*{.25cm}}

\newcounter{minutes}
\divide\time by 60
\newcounter{hours}
\multiply\time by 60
\addtocounter{minutes}{-\time}

\newcommand{\N}{\mathbf{N}} 
\newcommand{\C}{\mathbb{C}} 
\newcommand{\D}{\mathbb{D}} 
\newcommand{\R}{\mathbb{R}} 
\newcommand{\symD}{\Omega} 

\begin{document}

\title[Harmonic Shears of Slit and Polygonal Mappings]{Harmonic Shears of Slit and Polygonal Mappings}
\date{\today}

\author[S.~Ponnusamy]{Saminathan Ponnusamy}
\email{samy@iitm.ac.in}
\address{Indian Institute of Technology Madras,
         Department of Mathematics,
         Chennai-600 036,
         India.}

\author[T.~Quach]{Tri Quach${}^\textrm{{\tiny\FiveFlowerOpen}}$}
\email{tri.quach@tkk.fi}
\address{Aalto University,
         Department of Mathematics and Systems Analysis,
         P.O. Box 11100,
         FI-00076 Aalto,
         Finland.}

\author[A.~Rasila]{Antti Rasila}
\email{antti.rasila@iki.fi}
\address{Aalto University,
         Department of Mathematics and Systems Analysis,
         P.O. Box 11100,
         FI-00076 Aalto,
         Finland.}


\thanks{This research was supported by a grant from the Jenny and Antti Wihuri Foundation. \\
${}^\textrm{{\tiny\FiveFlowerOpen}}$ The author was supported by a grant (ma2011n25) from the Magnus Ehrnrooth Foundation.}

\begin{abstract}
In this paper, we study harmonic mappings by using the {\it shear construction}, introduced by
Clunie and Sheil-Small in 1984. We consider two classes of conformal mappings, each of which
maps the unit disk $\D$ univalently onto a domain which is convex in the horizontal direction,
and shear these mappings with suitable dilatations $\omega$.
Mappings of the first class map the unit disk $\D$ onto  four-slit domains
and mappings of the second class take $\D$ onto regular $n$-gons.
In addition, we discuss the minimal surfaces associated with such harmonic mappings. Furthermore,
illustrations of mappings and associated minimal surfaces are given by using {\sc Mathematica}.
\end{abstract}

\maketitle




\section{Introduction}
A complex-valued harmonic function $f$ defined on the unit disk $\D$ is called a
{\it harmonic mapping} if it maps $\D$ univalently onto a domain
$\symD \subset \C$. Note that it is not required that the real and the imaginary part of $f$
satisfy the Cauchy-Riemann equations. 
In 1984, Clunie and Sheil-Small \cite{clunie} showed that many classical results for conformal mappings have natural analogues for harmonic mappings, and hence they can be regarded as a generalization of conformal mappings. 
Each harmonic mapping in $\D$
has a canonical presentation $f=h+\overline{g}$, where $h$ and $g$ are analytic in $\D$ and
$g(0)=0$. A harmonic mapping $f=h+\overline{g}$ is called {\it sense-preserving}
if the Jacobian $J_f = |h'|^2-|g'|^2$ is positive in $\D$. Then $f$ has an {\it analytic dilatation}
$\omega =g'/h'$ such that $|\omega(z)| < 1$ for $z\in\D$. For basic properties of harmonic
mappings we refer to \cite{Dur,pon-ras}.

\lisays
A domain $\symD \subset \C$ is said to be {\it convex in the horizontal direction} (CHD)
if its intersection with each horizontal line is connected (or empty).
A univalent harmonic mapping is called a CHD mapping if its range is a CHD domain.
Construction of a harmonic mapping $f$ with prescribed dilatation
$\omega$ can be done effectively by the {\it shear construction} the devised by
Clunie and Sheil-Small \cite{clunie}.

\begin{theorem}
Let $f = h+\overline{g}$ be a harmonic and locally univalent in the unit disk $\D$.
Then $f$ is univalent in $\D$ and its range is a {\rm CHD} domain if and only if $h-g$ is a conformal
mapping of $\D$ onto a {\rm CHD} domain.
\end{theorem}

Suppose that $\varphi$ is a CHD conformal mapping. For a given dilatation $\omega$, the harmonic shear
$f=h+\overline{g}$ of $\varphi$ is obtained by solving the differential equations
\[
\left\{ \begin{split}
h'-g' &=\varphi', \\
\omega h' - g' & = 0.
\end{split} \right.
\]
From the above equations, we obtain
\begin{equation} \label{eqn: fun-h}
h(z) = \int_0^z \frac{\varphi'(\zeta)}{1-\omega(\zeta)} \, d\zeta.
\end{equation}
For the anti-analytic part $g$, we have
\begin{equation} \label{eqn: fun-g}
g(z) = \int_0^z \omega(\zeta) \frac{\varphi'(\zeta)}{1-\omega(\zeta)} \, d\zeta.
\end{equation}
Observe that
\begin{equation} \label{eqn: fun-f}
f(z) = h(z) +\overline{g(z)}
      = 2 \, \textrm{Re} \left[ \int_0^z \frac{\varphi'(\zeta)}{1-\omega(\zeta)} \, d\zeta \right] - \overline{\varphi(z)}.
\end{equation}
We shall use (\ref{eqn: fun-h}) to find the analytic part $h$ of the harmonic mapping
$f=h+\overline{g}$. Then the anti-analytic part $g$ of the harmonic mapping $f$ can be obtained from the identity $g = h-\varphi$, or computed via (\ref{eqn: fun-f}).

\lisays
It is known that the class of harmonic mappings has a close connection with the theory of minimal surfaces.
In the space $\R^3$, the {\it minimal surface} is a surface which minimizes the area with a fixed curve
as its boundary. This minimization problem is called {\it Plateau's Problem}. Discussion concerning the differential
geometric approach to the subject can be found from the book by Pressley \cite{pressley}.

\lisays
Our results concerning minimal surfaces are based on the Weierstrass-Enneper representation. Let $S$ be a non-parametric minimal surface over a simply connected domain $\symD$ in $\C$ given by
\[
S=\{(u, v, F(u,v) ):\, u+iv \in \symD\},
\]
where $(u,v)$ identifies the complex plane $\R^2$, which underlies the domain of $F$.
The following result is known as the Weierstrass-Enneper representation. This representation provides a link between
harmonic univalent mappings and minimal surfaces. The surface $S$ is minimal surfaces
if and only if $S$ has the following representation
\[
S=\left\{ \left(\textrm{Re} \, \int_0^z \varphi_1(\zeta)\,d\zeta+c_1, \,
\textrm{Re}\, \int_0^z \varphi_2(\zeta)\,d\zeta+c_2,\, \textrm{Re} \, \int_0^z \varphi_3(\zeta)\,d\zeta+c_3 \right):\, z\in\D\right\},
\]
where $\varphi_1,\, \varphi_2,\, \varphi_3$ are analytic such that $\varphi_1^2+\varphi_2^2+\varphi_3^2 =0$,
and
\[
f(z)=u(z)+iv(z)  =\textrm{Re} \, \int_0^z \varphi_1(z) \, dz+i \textrm{Re} \, \int_0^z \varphi_2(z) \, dz+c
\]
is a sense-preserving univalent harmonic mapping from $\D$ onto $\symD$. In this case,
the surface $S$ is called a minimal surface over $\symD$ with the projection $f=u+iv$. Further
information about  the relation between harmonic mappings and minimal surfaces can be found from the book of Duren \cite{Dur}.

\lisays
Systematical construction of harmonic shears of mappings of the unit disk and unbounded strip domains, and their boundary behaviour are presented in the article by Greiner \cite{greiner}. In most cases the dilatation is chosen to be $\omega(z) = z^n$.

\lisays
In this paper, we study two classes of conformal mappings, each of which map $\D$ univalently onto
a domain which is convex in the horizontal direction. The first one involves the mapping
\[
\varphi(z) = A \log\left( \frac{1+z}{1-z} \right) + B \frac{z}{1+cz+z^2},
\]
which maps $\D$ onto $\C$ minus four symmetric half-lines. In \cite{gan-wid}, Ganczar and Widomski have
studied some special cases of this mapping and its harmonic shears. Analytic examples of harmonic
shears of $\varphi$ with dilatations
\[
\omega (z) = \frac{1-z^{2k}}{1+z^{2k}},  \quad k=1,2,
\]
along with illustrations, are given in \cite{dnw}.

The second case is related to the conformal mapping (see \cite[p. 196]{nehari})
\begin{equation} \label{eqn:n-gonmap}
\varphi(z) = \int_0^z (1-\zeta^n)^{-2/n} \, d\zeta,
\end{equation}
which maps the unit disk $\D$ onto a regular $n$-gon. In \cite{driver}, Driver and Duren discussed the harmonic shear of $\varphi$ by the choice of the dilatations
$\omega (z) = z^n$, $\omega(z)=z$ and $\omega(z)=z^{n/2}$. With the last choice of the dilatation,
it is assumed that $n$ is even. The minimal surfaces of these harmonic shears were studied as well.

\lisays
The outline of the paper is as follows: In the first part, we study the four-slit conformal
mapping, and its harmonic shears, with the dilatation $\omega(z) = z^n$. Then we consider
conformal mappings of regular $n$-gons, with dilatations $\omega(z) = z^{2n}$ and $\omega(z) = z^2$.
These dilatations lead naturally to the Appel's hypergeometric functions, which are formal
extensions of the Gaussian hypergeometric functions into two variables. Finally, we
consider minimal surfaces obtained by shearing conformal mapping of
the regular $n$-gon described by \eqref{eqn:n-gonmap}. The results are also illustrated by using {\sc Mathematica}.


\section{Shearing of Four Slit Conformal Mappings}
In this section we shall give examples of harmonic shear of unbounded conformal mappings with
a suitable dilatations. For $A,B>0$ and $c \in [-2,2]$, let us consider the function $\varphi$ defined by
\begin{equation} \label{eqn: fun-slit}
\varphi(z) = A \log\left( \frac{1+z}{1-z} \right) + B \frac{z}{1+cz+z^2}.
\end{equation}
Recall from \cite{dnw} that $\varphi$ is univalent and it  maps the unit disk $\D$ onto a
domain convex in the direction of the real axis. In special cases, namely $c=-2$ and $c=2$, the
image of the unit disk $\D$ under $\varphi$ is
\[
\C \backslash \left\{ x \pm \frac{A\pi}{2}i : x \in \left( - \infty, \frac{A}{2} \log \frac{2A}{B} - \frac{2A+B}{4} \right]\right\},
\]
and
\[
\C \backslash \left\{ x \pm \frac{A\pi}{2}i : x \in  \left[ -\frac{A}{2} \log \frac{2A}{B} + \frac{2A+B}{4}, \infty  \right)\right\},
\]
respectively.
In the case $c=0$, the mapping $\varphi$ maps the unit disk $\D$ onto $\C$ minus the following half-lines:
\[
\left\{ x \pm \frac{A\pi}{2}i : x \in  \left( - \infty, -\frac{A}{2} \log \frac{\sqrt{2A+B}+B}{\sqrt{2A+B}-B} - \frac{\sqrt{B(2A+B)}}{4} \right]\right\},
\]
and
\[
\left\{ x \pm \frac{A\pi}{2}i : x \in  \left[ \frac{A}{2} \log \frac{\sqrt{2A+B}+B}{\sqrt{2A+B}-B} + \frac{\sqrt{B(2A+B)}}{4}, \infty  \right)\right\}.
\]
By writing $c = -2\cos \gamma$, $\gamma \in (0,\pi)$, the equation (\ref{eqn: fun-slit}) takes the form
\begin{equation} \label{eqn: fun-slit-2}
\varphi(z) = A \log\left( \frac{1+z}{1-z} \right) + B \frac{z}{(1- e^{i\gamma}z)(1- e^{-i\gamma}z)}.
\end{equation}
Also in \cite{dnw}, the authors considered harmonic shears of $\varphi$ for some choices of $A$ and $B$ with dilatations
\[
\omega(z) = \frac{1-z^{2k}}{1+z^{2k}}, \quad k=1,2.
\]
Let us now consider the harmonic shear of $\varphi$ defined in (\ref{eqn: fun-slit-2}) with the dilatation
$\omega(z) = z^n$, $n \geq 2$. First, we need the derivative of $\varphi$ and the direct calculation gives
\begin{equation} \label{eqn: fun-slit-der}
\varphi'(z) = \frac{2A}{1-z^2} - \frac{B}{2\sin \gamma}i \left[ \frac{e^{i\gamma}}{(1-e^{i\gamma}z)^2} - \frac{e^{-i\gamma}}{(1-e^{-i\gamma}z)^2}\right].
\end{equation}
Thus by (\ref{eqn: fun-h}) and (\ref{eqn: fun-slit-der}), we have
\[ \begin{split}
h(z) = &\,\, 2A \int_0^z \frac{d\zeta}{(1-\zeta^2)(1-\zeta^n)}  \\
      &\,\, - \frac{B}{2\sin \gamma}i \left[ \int_0^z \frac{e^{-i\gamma}\, d\zeta}{(\zeta-e^{-i\gamma})^2 (1-\zeta^n)}
      - \int_0^z \frac{e^{i\gamma}\, d\zeta}{(\zeta-e^{i\gamma})^2(1-\zeta^n)}\right].
\end{split}
\]
We shall write this in the form
\[
h(z) = 2AI_1 - \frac{B}{2\sin \gamma}i \left( e^{-i\gamma} I_2 - e^{i\gamma} I_3 \right),
\]
where
\begin{align*}
I_1 & =\int_0^z \frac{d\zeta}{(1-\zeta^2)(1-\zeta^n)},\\
I_2 & =\int_0^z \frac{d\zeta}{(\zeta-e^{-i\gamma})^2 (1-\zeta^n)},\\
I_3 &=\int_0^z \frac{d\zeta}{(\zeta-e^{i\gamma})^2(1-\zeta^n)}.
\end{align*}
By partial fractions and $z_k = e^{2\pi ik/n}$, $k=0, \cdots, n-1$, we have
\[
\frac{1}{1-z^n} = -\frac{1}{n}\sum_{k=0}^{n-1} \frac{z_k}{z - z_k}.
\]
Therefore the first integral $I_1$ can be rewritten as
\[
I_1 = \frac{1}{n}I_{1,0} + \frac{1}{n}\sum_{k=1}^{n-1} I_{1,k} = \frac{1}{n} \int_0^z \frac{d\zeta}{(1-\zeta)^2(1+\zeta)} +
 \frac{1}{n} \sum_{k=1}^{n-1} \int_0^z \frac{z_k \, d\zeta}{(1-\zeta)(1+\zeta)(z_k-\zeta)}.
\]
We remark that for the case $n=1$, the latter part of the integral inside the summation sign should be omitted.
For the first integral, we get
\[
I_{1,0} = \int_0^z \frac{d\zeta}{(1-\zeta)^2(1+\zeta)} = \frac{1}{2}\frac{z}{1-z} + \frac{1}{4} \log \left( \frac{1+z}{1-z} \right).
\]
The latter integral $I_{1,k}$ depends on whether $n$ is odd or even. Assuming that $n$ is odd, we easily see that
\begin{align*}
I_{1,k} & = \int_0^z \frac{z_k \, d\zeta}{(1-\zeta)(1+\zeta)(z_k-\zeta)} \\
     & =\frac{z_k}{2(z_k-1)} \int_0^z \frac{d\zeta}{1-\zeta} + \frac{z_k}{2(z_k+1)} \int_0^z \frac{d\zeta}{1+\zeta} - \frac{z_k}{z_k^2-1} \int_0^z \frac{d\zeta}{z_k-\zeta}.
\end{align*}
Thus
\[
I_{1,k}= -\frac{z_k \log(1-z)}{2(z_k-1)}  + \frac{z_k \log(1+z)}{2(z_k+1)}  + \frac{z_k}{z_k^2-1} \log \left( \frac{z_k-z}{z_k}\right).
\]
Note that, by assumption, $z_k\neq \pm 1$.
To simplify our notation, let $\N = \{0, 1, \cdots, n-1 \}$ be an index set. Suppose that $a \in \N $, then we define $\N_a = \N \backslash \{a\}$, and $\N_{a,b} = \N_a \cap \N_b$.  With this notation, in case $n$ is even, we have
\[
\sum_{k=1}^{n-1} I_{1,k}  = \sum_{k\in\N_0} I_{1,k} = \sum_{k\in \N_{0,n/2}} I_{1,k} + I_{1,n/2},
\]
where
\[
I_{1,n/2} = \int_0^z \frac{d\zeta}{(1+\zeta)^2 (1-\zeta)} = \frac{1}{2}\frac{z}{1+z} +\frac{1}{4} \log \left( \frac{1+z}{1-z} \right).
\]
Next we compute the integrals $I_2$ and $I_3$. Assuming that $\eta = e^{i\gamma} \not= z_k$, i.e, $\gamma \not= 2\pi k/n$, $k=0, \cdots, n-1$, we compute
\begin{align*}
I_\eta & =\int_0^z \frac{d\zeta}{(\zeta-\eta)^2(1-\zeta^n)} \\
    & = -\frac{1}{n} \sum_{k=0}^{n-1} \int_0^z \frac{z_k \, d\zeta}{(\zeta-\eta)^2(\zeta-z_k)} \\
    & = -\frac{1}{n} \sum_{k=0}^{n-1} \left[ \frac{1}{\eta-z_k} \int_0^z \frac{d\zeta}{(\zeta-\eta)^2} + \frac{1}{(\eta-z_k)^2} \int_0^z \frac{d\zeta}{\eta-\zeta}  - \frac{1}{(\eta-z_k)^2} \int_0^z \frac{d\zeta}{z_k-\zeta} \right] \\
    & = -\frac{1}{n} \sum_{k=0}^{n-1} \left\{ \frac{1}{\eta-z_k} \left( \frac{1}{\eta-z} - \frac{1}{\eta} \right)  - \frac{1}{(\eta-z_k)^2} \left[ \log \left( \frac{\eta-z}{\eta} \right) - \log \left( \frac{z_k-z}{z_k} \right) \right] \right\}.
\end{align*}
In the case of $\gamma = 2\pi m/n$, for $m=0,\cdots, n-1$, we have
\begin{align*}
I_{3,m} & = \int_0^z \frac{d\zeta}{(\zeta-z_m)^2(1-\zeta^n)}  \\
    & = \frac{1}{n} \sum_{k\in \N_m} \int_0^z \frac{d\zeta}{(\zeta-z_m)^2(\zeta-z_k)} + \frac{1}{n}\int_0^z \frac{d\zeta}{(\zeta-z_m)^3}.
\end{align*}
The sum can be computed as above and the last integral is
\[
\frac{1}{n}\int_0^z \frac{d\zeta}{(\zeta-z_m)^3} = \frac{1}{2n} \left[ \frac{1}{(z-z_m)^2} -\frac{1}{z_m^2} \right].
\]
Note that, we have an identity $I_2 = I_{3,n-m}$ and if $m=0$, we have $I_2 = I_{3,0}$.

\lisays
Finally, the function $g$ can be solved readily from the following identity
\[
g = h - \varphi.
\]
In Table \ref{table: slit-map-I1} we have the integral $I_1$ of $h$, which depends on whether $n$ is odd or even. On the other hand, integrals $I_2$ and $I_3$ depend on whether $\gamma \not= 2\pi m/n$ or $\gamma = 2\pi m/n$, for an $m=0,\cdots, n-1$ and the result is given in Table \ref{table: slit-map-I2}.
\begin{table}[ht]
\caption{Analytic part $h = 2AI_1 + \frac{B}{2 \sin \gamma}i (e^{-i\gamma} I_2 - e^{i\gamma} I_3)$ of the harmonic shear $f$ of the slit mapping $\varphi$ with a dilatation $\omega(z)=z^n$. Note that in the case of $n=1$, the summation should be omitted. See Table \ref{table: slit-map-I2} for the integrals $I_2$ and $I_3$.}\label{table: slit-map-I1}.

\begin{tabular}{|c|c|}
\hline
$n$ & $I_1$  \\ \hline
even &  $\displaystyle\frac{1}{n} \left( I_{1,0} + \displaystyle\sum_{k=1}^{n-1} I_{1,k} \right)$ \\ \hline
odd  &  $\displaystyle\frac{1}{n} \left( I_{1,0} + I_{1,n/2} + \displaystyle\sum_{k\in\N_{0,n/2}} I_{1,k} \right)$ \\
\hline
\end{tabular}
\end{table}

\begin{table}[ht]
\caption{The integrals $I_2$ and $I_3$ for the analytic part $h = 2AI_1 + \frac{B}{2 \sin \gamma}i (e^{-i\gamma} I_2 - e^{i\gamma} I_3)$ of the harmonic shear $f$ of the slit mapping $\varphi$ with a dilatation $\omega(z)=z^n$. See Table \ref{table: slit-map-I1} for the integral $I_1$.}\label{table: slit-map-I2}
\begin{tabular}{|c|c|c|}
\hline
$\gamma$ & $I_2$ & $I_3$ \\ \hline
is not $2\pi m/n$  & $I_\eta$ & $I_{\overline{\eta}}$ \\ \hline
is $2\pi m/n$ & $I_{3,n-m}$ & $I_{3,m}$ \\
\hline
\end{tabular}
\end{table}

\lisays
Therefore, we have obtained the following result:

\begin{theorem}
Let $\varphi$ be given by \eqref{eqn: fun-slit-2}. Then the harmonic shear $f = h + \overline{g}$,
where $h$ is given in Table \ref{table: slit-map-I1} and Table \ref{table: slit-map-I2} and the anti-analytic
part $g = h-\varphi$, maps the unit disk $\D$ univalently onto a domain which is convex in the horizontal direction.
\end{theorem}

\lisays
In the case of $A = \frac{1}{2} \sin^2\alpha$, $B = \cos^2\alpha$ and $c=0$, we have the Corollary.
\begin{corollary} \label{cor: nehari}
Let $\varphi_\alpha$ be defined as
\[
\varphi_\alpha(z) = \frac{1}{2} \sin^2\alpha \log\left( \frac{1+z}{1-z}\right) + \frac{z \cos^2\alpha}{1+z^2}.
\]
Then $\varphi_\alpha$ maps the unit disk $\D$ univalently onto $\C$ minus the following half-lines
\[
\left\{ x \pm \frac{\pi\sin^2\alpha}{4}i : x \in  \left( - \infty, -\frac{1}{2}\sin^2\alpha \log \cot \frac{\alpha}{2} - \frac{\cos\alpha}{2} \right]\right\},
\]
and
\[
\left\{ x \pm \frac{\pi\sin^2\alpha}{4}i : x \in  \left[ \frac{1}{2}\sin^2\alpha \log \cot \frac{\alpha}{2} + \frac{\cos\alpha}{2}, \infty  \right)\right\},
\]
which is convex in the horizontal direction. Then the harmonic shear $f_\alpha = h + \overline{g}$, with
\[
h = \sin^2\alpha \,I_1 - \frac{\cos^2\alpha}{2} (I_2 - I_3),
\]
where $I_1, I_2,I_3$ are given in Table \ref{table: slit-map-I1} and Table \ref{table: slit-map-I2}, and
$g = h-\varphi_a$, maps the unit disk $\D$ onto a {\rm CHD} domain.
\end{corollary}

\begin{figure}[!ht]
\begin{center}
\includegraphics[width=0.5\textwidth]{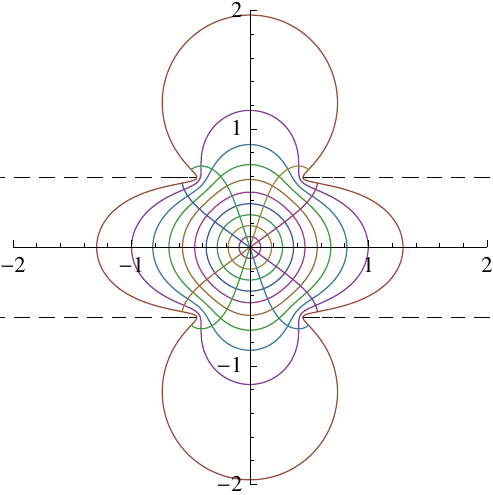}
\caption{Conformal mapping $\varphi_\alpha$ of the unit disk $\D$ onto a four-slit domain with $\alpha= \pi/3$, and the omitted half-lines are dashed. Illustrations of harmonic shears are given in Figure \ref{fig: nehari-pi3}.
} \label{fig: conf-nehari}
\end{center}
\end{figure}

\begin{figure}[!ht]
\begin{center}
\subfloat[Dilatation $\omega(z) = z$.]{\parbox{.45\textwidth}{\centering\includegraphics[width=.45\textwidth]{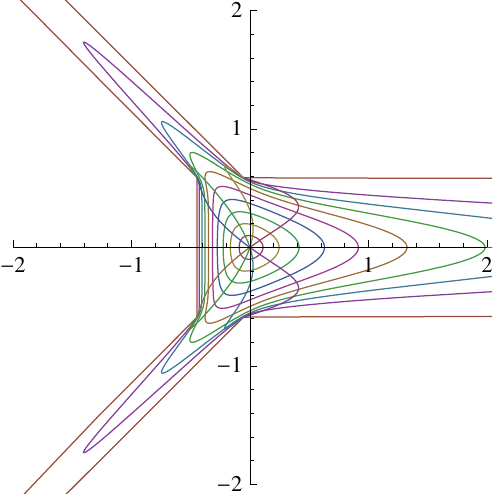}}}
\hspace*{.5cm}
\subfloat[Dilatation $\omega(z) = z^2$.]{\parbox{.45\textwidth}{\centering\includegraphics[width=.45\textwidth]{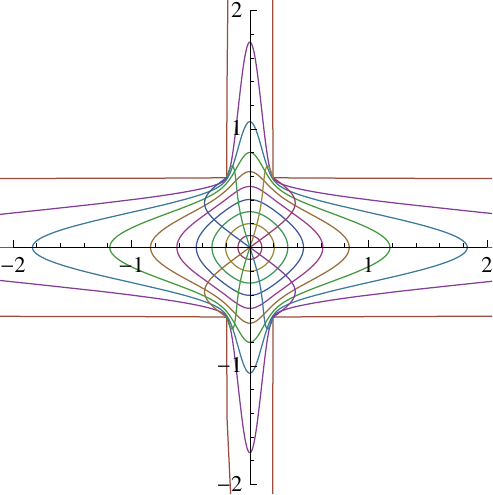}}}
\caption{Harmonic shears $f_\alpha$ given in Corollary \ref{cor: nehari} with $\alpha= \pi/3$ and a dilatation $\omega(z)=z^n$, for $n=1,2$.
} \label{fig: nehari-pi3}
\end{center}
\end{figure}

\lisays
The above conformal mapping $\varphi_\alpha$ is given in \cite[p. 197]{nehari}.
In Figure \ref{fig: conf-nehari}, we have the conformal mapping $\varphi_\alpha$ of the unit disk $\D$ onto
a four-slit domain with $\alpha=\pi/3$. Harmonic shears of this mapping as given in
Corollary \ref{cor: nehari} are shown in Figure \ref{fig: nehari-pi3} with $\alpha= \pi/3$, for $n=1,2$.

\section{Shearing of Regular $n$-gons} \label{section: polygon}
As pointed out in the introduction, the authors in \cite{driver}
considered the conformal mapping $\varphi$ given by \eqref{eqn:n-gonmap}.
Then the harmonic shear $f$ (with dilatations $\omega(z) = z^{kn/2}$, $k=1,2$, and $\omega(z)=z$)
can be given in terms of the Gaussian hypergeometric function, which is defined by
\[
F(a,b;c;z) = 1 + \sum_{n=1}^\infty \frac{(a)_n (b)_n}{n! (c)_n} z^n, \quad |z| < 1,
\]
where
\[
(\alpha)_n = \alpha (\alpha+1) \cdots (\alpha+n-1) = \frac{\Gamma(\alpha+n)}{\Gamma(\alpha)}, \quad \alpha \in \C,
\]
is the Pochhammer symbol. For $\textrm{Re}\, c > \textrm{Re}\, b > 0$, this can also be written as the Euler integral
\[
F(a,b,c;z) = \frac{\Gamma(c)}{\Gamma(b)\Gamma(c-b)} \int_0^1 t^{b-1} (1-t)^{c-b-1} (1-zt)^{-a} \, dt.
\]
Now we consider the dilatation $\omega(z) = z^{2n}$. Therefore, by (\ref{eqn: fun-h}),
(\ref{eqn: fun-g}) and (\ref{eqn:n-gonmap}), we have
\begin{equation} \label{eqn: fun-hg-w2} \left\{ \begin{split}
h(z) & = \int_0^z (1-\zeta^n)^{-2/n} (1-\zeta^{2n})^{-1} \, d\zeta, \\
g(z) & = \int_0^z \zeta^{2n}(1-\zeta^n)^{-2/n} (1-\zeta^{2n})^{-1} \, d\zeta.
\end{split} \right.
\end{equation}
Notice that we can write $h$ and $g$ as
\[ \left\{ \begin{split}
h(z) & = \int_0^z (1-\zeta^n)^{-1-2/n} (1+\zeta^n)^{-1} \, d\zeta, \\
g(z) & = \int_0^z \zeta^{2n}(1-\zeta^n)^{-1-2/n} (1+\zeta^n)^{-1} \, d\zeta.
\end{split} \right.
\]
Using the change of variable $\zeta = t^{1/n}z$, we obtain
\[
h(z) = \frac{z}{n} \int_0^{1} (1-z^nt)^{-1-2/n} (1+z^nt)^{-1} t^{1/n-1} \, dt.
\]
In order to rewrite $h$ in a compact form for our purpose, we recall the first Appel's hypergeometric
function \cite[p. 73]{bailey}, which is defined by
\[
F_1(a,b_1,b_2;c;x,y) = \sum_{k=0}^\infty \sum_{l=0}^\infty \frac{(a)_{k+l} (b_1)_k (b_2)_l}{(c)_{k+l} \,k!\, l!} x^k y^l,
\]
where $(\alpha)_n$ is the Pochhammer symbol given above. As for hypergeometric functions,
Appel's hypergeometric functions can be defined by Euler's integral as follows \cite[p. 77]{bailey}:
\[
F_1(a,b_1,b_2;c;x,y) = \frac{\Gamma(c)}{\Gamma(a)\Gamma(c-a)} \int_0^1 t^{a-1} (1-t)^{c-a-1} (1-xt)^{-b_1} (1-yt)^{-b_2} \, dt,
\]
where $\textrm{Re}\, c > \textrm{Re}\, a > 0$. Set
$a = \frac{1}{n}, b_1 = 1+\frac{2}{n}, b_2 = 1, c = 1+\frac{1}{n}$. We have
\begin{equation} \label{eqn: polygon-h-z2n}
h(z) = z F_1\left( \frac{1}{n}, 1+\frac{2}{n}, 1; 1+\frac{1}{n}; z^n, -z^n \right).
\end{equation}
A direct computation shows that $g$ defined by (\ref{eqn: fun-hg-w2}) can be written in the form
\begin{equation} \label{eqn: polygon-g-z2n}
g(z) = \frac{z^{2n+1}}{2n+1} F_1\left( 2 +\frac{1}{n}, 1 +\frac{2}{n}, 1; 3 +\frac{1}{n}; z^n, -z^n \right).
\end{equation}
This proves the following result.

\begin{theorem} \label{thm: polygon-2n}
Let $\varphi$ be given by \eqref{eqn:n-gonmap}. Then the harmonic shear $f = h + \overline{g}$,
where $h$ and $g$ are given in \eqref{eqn: polygon-h-z2n} and \eqref{eqn: polygon-g-z2n},
respectively, maps the unit disk $\D$ univalently onto a domain which is convex in the horizontal direction.
\end{theorem}

\begin{figure}[!hb]
\begin{center}
\includegraphics[width=0.43\textwidth]{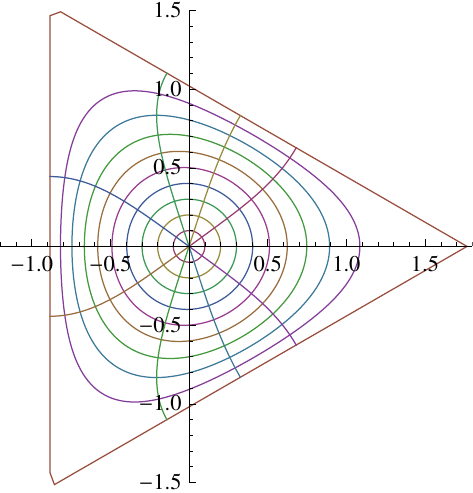}
\includegraphics[width=0.43\textwidth]{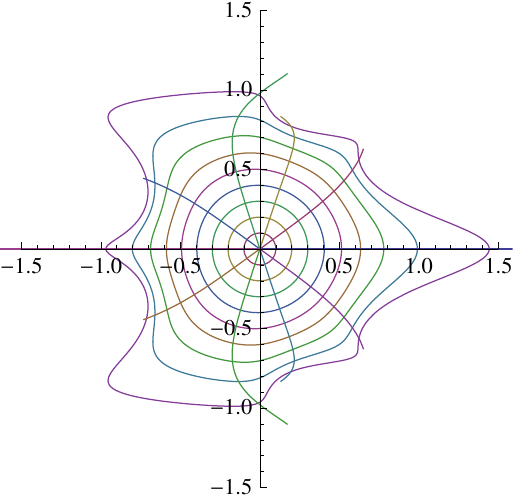}
\includegraphics[width=0.43\textwidth]{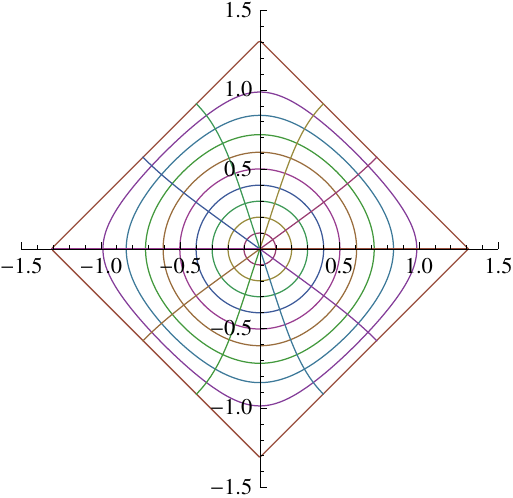}
\includegraphics[width=0.43\textwidth]{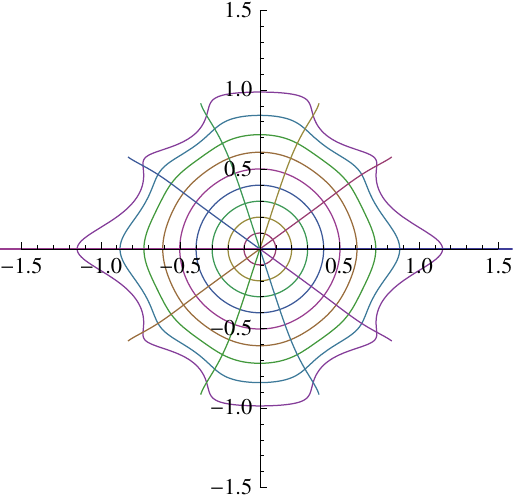}
\includegraphics[width=0.43\textwidth]{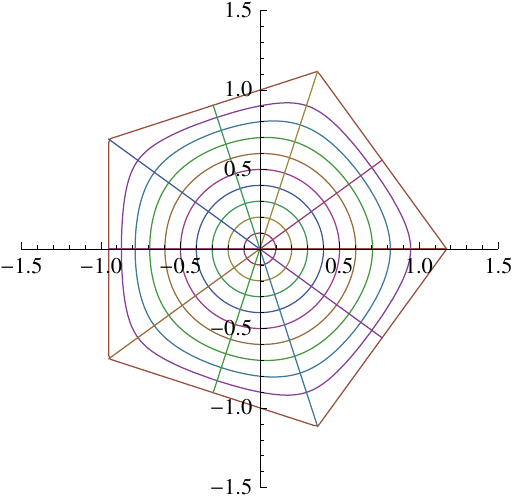}
\includegraphics[width=0.43\textwidth]{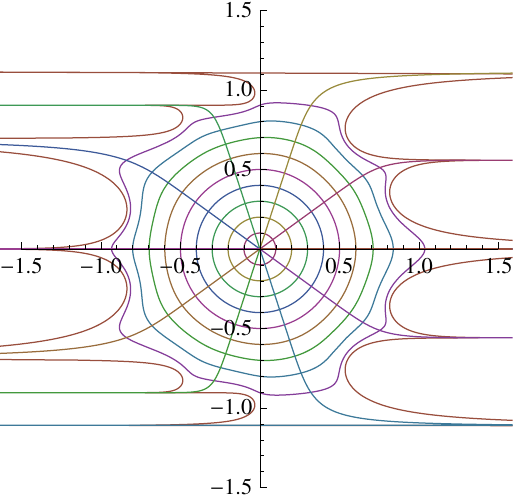}
\caption{Conformal mappings from the unit disk $\D$ onto a regular $n$-gon and its harmonic shears with the dilatation $\omega(z) = z^{2n}$, $n=3,4,5$ For a comparison with a dilatation $\omega(z)=z^{n/2}$, see \cite{driver}.
} \label{fig: polygon-double}
\end{center}
\end{figure}

\lisays
In Figure \ref{fig: polygon-double}, we have illustrations of the conformal mapping $\varphi$
onto regular $n$-gon and the harmonic shear $f = h +\overline{g}$ with dilatation $\omega(z) = z^{2n}$, for $n=3,4,5$.
It is worth remarking that this situation was not considered by Duren and Driver in \cite{driver}.

\begin{figure}[!Ht]
\centering
\subfloat[$n=3$]{\parbox{.8\textwidth}{\centering\includegraphics[width=.55\textwidth]{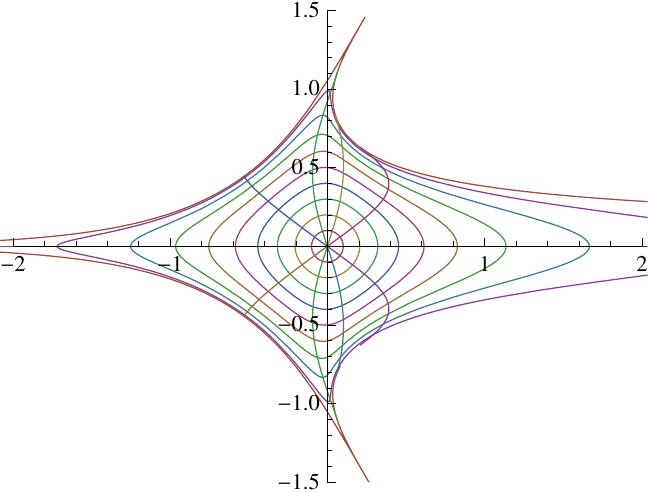}}}\\
\subfloat[$n=5$]{\parbox{.8\textwidth}{\centering\includegraphics[width=.55\textwidth]{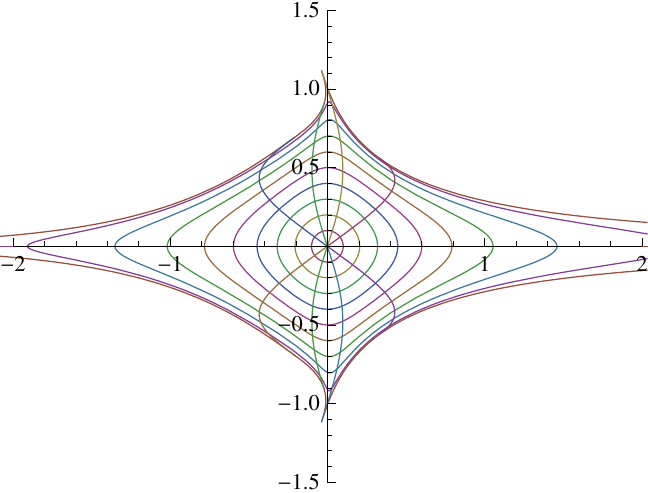}}}
\caption{Harmonic shear of the mapping $\varphi$, which maps the unit disk $\D$ onto a regular $n$-gon with the dilatation $\omega(z)=z^2$.} \label{fig: polygon-square}
\end{figure}

\clearpage
Next we let $n$ be odd and we consider the dilatation $\omega(z) = z^2$. Thus, by (\ref{eqn: fun-h}), (\ref{eqn: fun-g}) and (\ref{eqn:n-gonmap}), we have
\[ \left\{ \begin{split}
h(z) & = \int_0^z (1-\zeta^n)^{-2/n} (1-\zeta^2)^{-1} \, d\zeta, \\
g(z) & = \int_0^z \zeta^2 (1-\zeta^n)^{-2/n} (1-\zeta^2)^{-1} \, d\zeta.
\end{split} \right.
\]
Since $n$ is assumed to be odd, we have
\begin{align*}
\frac{1+\zeta^n}{1+\zeta} \frac{1-\zeta^n}{1-\zeta}
&= (1 - \zeta + \cdots - \zeta^{n-2} + \zeta^{n-1} )(1 + \zeta + \cdots + \zeta^{n-1} )\\
& = 1 + \zeta^2 + \cdots + \zeta^{2(n-1)},
\end{align*}
and as a consequence of this observation, $h$ defined above takes the form
\[
h(z) = \int_0^z \frac{1 + \zeta^2 + \cdots + \zeta^{2(n-1)}}{(1-\zeta^n)^{1+2/n} (1+\zeta^n)} \, d\zeta = \sum_{k=0}^{n-1} \int_0^z \frac{\zeta^{2k}}{(1-\zeta^n)^{1+2/n} (1+\zeta^n)} \, d\zeta.
\]
A similar expression holds for $g$ as well. Finally, by computation, we obtain that
\begin{equation} \label{eqn: polygon-h-g}
\left\{ \begin{split}
h(z) & = \sum_{k=0}^{n-1} \frac{z^{2k+1}}{2k+1} F_1\left( \frac{2k+1}{n}, 1+\frac{2}{n}, 1; 1+ \frac{2k+1}{n}, z^n,-z^n \right), \\
g(z) & = \sum_{k=1}^{n} \frac{z^{2k+1}}{2k+1} F_1\left( \frac{2k+1}{n}, 1+\frac{2}{n}, 1; 1+ \frac{2k+1}{n}, z^n,-z^n \right).
\end{split} \right.
\end{equation}
We have shown the following result:

\begin{theorem} \label{thm: polygon-z2}
Let $\varphi$ be given by \eqref{eqn:n-gonmap}. Then the harmonic shear $f = h + \overline{g}$, where $h$ and $g$ are given in \eqref{eqn: polygon-h-g}, maps the unit disk $\D$ univalently onto a domain which is convex in the horizontal direction.
\end{theorem}
In Figure \ref{fig: polygon-square}, illustrations of the conformal mapping $\varphi$ onto a regular $n$-gon and the harmonic shear $f$ with dilatation $\omega(z) = z^2$, for $n=3,5$, are given.


\section{Minimal Surfaces}

\begin{figure}[!bt]
\begin{center}
\subfloat[$n=3$]{\centering\includegraphics[width=.65\textwidth]{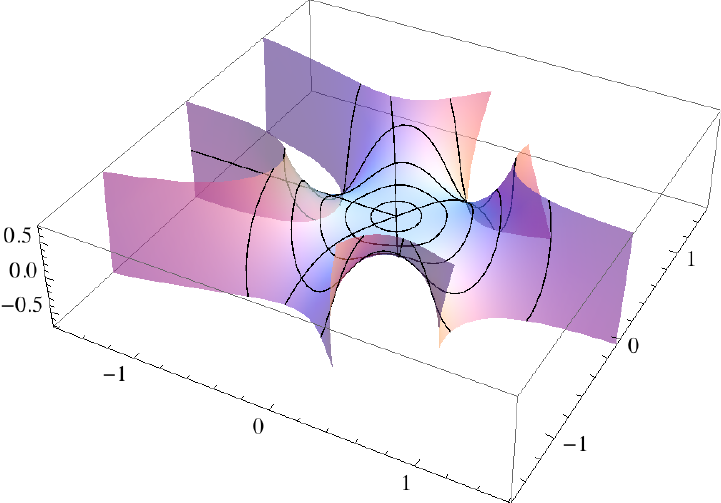}} \\
\subfloat[$n=4$]{\centering\includegraphics[width=.65\textwidth]{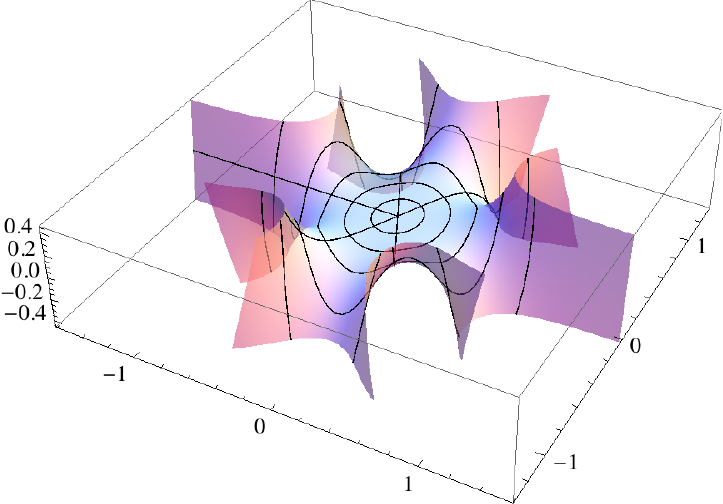}}
\caption{Minimal surfaces of the harmonic shear $f$, which is given by Theorem \ref{thm: polygon-2n}.
} \label{fig: minimal-2n}
\end{center}
\end{figure}

A harmonic function $f = h + \overline{g}$ can be lifted to a minimal surface if and only if the
dilatation $\omega$ is the square of an analytic function. Suppose that $\omega = q^2$ for some analytic
function $q$ in the unit disk $\D$. Then the corresponding minimal surfaces has the form
\[
\{u,v,w\} = \{\textrm{Re}\, f,  \textrm{Im}\, f, 2\,\textrm{Im}\, \psi\},
\]
where
\[
\psi(z) = \int_0^z q(\zeta) \frac{\varphi(\zeta)}{1-\omega(\zeta)} \, d\zeta.
\]
Let $\varphi$ be a conformal mapping, which maps the unit disk $\D$ onto a regular $n$-gon and let the dilatation be $\omega(z) = z^{2n}$. Then the minimal surfaces is determined by the integral
\[
\psi(z) = \int_0^z \zeta^n (1-\zeta^n)^{-1-2/n} (1+\zeta^n)^{-1} \, d\zeta.
\]
By the substitution $\zeta= t^{1/n}z$, we have
\[
\psi(z) = \frac{z^{n+1}}{n} \int_0^1 t^{1/n} (1-z^n t)^{-1-2/n} (1+z^n t)^{-1} \, dt.
\]
Appel's hypergeometric presentation gives
\[
\psi(z) = z^{n+1} F_1 \left( 1+\frac{1}{n}, 1+\frac{2}{n}; 1; 2+\frac{1}{n}; z^n, -z^n \right).
\]
In Figure \ref{fig: minimal-2n}, we have the minimal surfaces for the above mapping for $n=3,4$.

\begin{figure}[H]
\centering
\subfloat[$n=3$]{\centering\includegraphics[width=.65\textwidth]{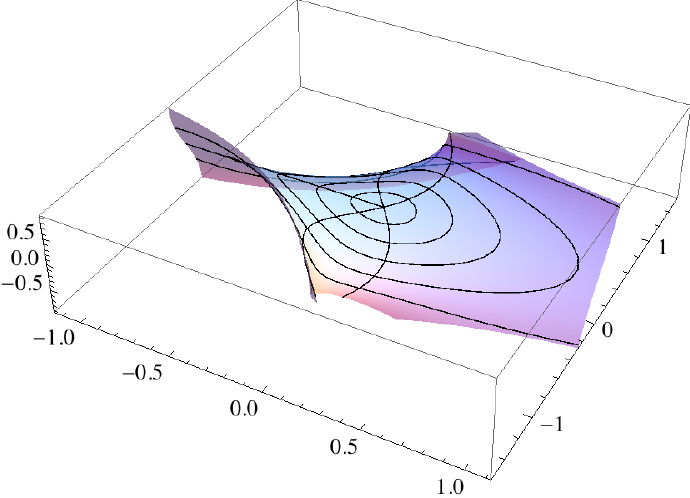}}\\
\subfloat[$n=5$]{\centering\includegraphics[width=.65\textwidth]{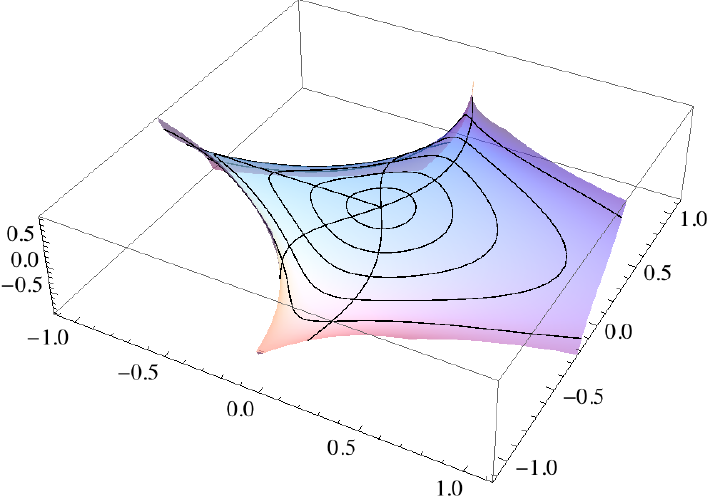}}
\caption{Minimal surfaces of the harmonic shear $f$, which is given by Theorem \ref{thm: polygon-z2}.
} \label{fig: minimal-z2}
\end{figure}

\lisays
In the second polygonal example, let $n$ be odd and let the dilatation be $\omega(z)=z^2$. In this case the minimal surface lifting is given by integral
\begin{align*}
\psi(z) & = \int_0^z \zeta (1-\zeta^n)^{-2/n} (1-\zeta^2)^{-1} \, d\zeta \\
    & = \sum_{k=0}^{n-1} \int_0^z \zeta^{2k+1} (1-\zeta^n)^{-1-2/n} (1+\zeta^n)^{-1} \, d\zeta \\
    & = \sum_{k=0}^{n-1} \frac{z^{2(k+1)}}{n}\int_0^1 t^{2(k+1)/n -1} (1-z^n t)^{-1-2/n} (1+z^n t)^{-1} \, dt.
\end{align*}
Again the above integral can be written by Appel's hypergeometric function
\begin{equation} \label{eqn: minimal-z2}
\psi(z) = \sum_{k=0}^{n-1} \frac{z^{2(k+1)}}{2(k+1)} F_1 \left( \frac{2(k+1)}{n}, 1 +\frac{2}{n}; 1; 1+\frac{2(k+1)}{n}; z^n, -z^n\right).
\end{equation}
In Figure \ref{fig: minimal-z2}, minimal surfaces of \eqref{eqn: minimal-z2}, for $n=3,5$, are illustrated.
\clearpage


%



\end{document}